\documentclass[12pt]{amsart}
 \usepackage{mathtools}
\mathtoolsset{showonlyrefs,showmanualtags} 

\usepackage{amssymb, amsmath, amsthm, esint}
\usepackage{mathrsfs}
\usepackage{bbm,dsfont}
\usepackage{hyperref}
\usepackage{mathtools}
\usepackage{fullpage}
\usepackage{color}

\usepackage[backrefs]{amsrefs}

 \usepackage[top=1.5in, bottom=1.5in, left=1.25in, right=1.25in]	{geometry}
\usepackage[all]{xy}
\usepackage{amsmath}
\usepackage{amssymb}
\usepackage{amsthm}
\usepackage{amscd}

\newcommand{\J}{\mathcal{J}}
\newcommand{\RR}{\mathbb  R}

\numberwithin{equation}{section}

\newtheorem{theorem}[equation]{Theorem}
\newtheorem{definition}[equation]{Definition}
\newtheorem{proposition}[equation]{Proposition}
\newtheorem{cor}[equation]{Corollary}
\newtheorem{lemma}[equation]{Lemma}

\DeclareMathOperator{\supp}{supp}

\usepackage{hyperref} 
\hypersetup{
    colorlinks=true,       
    linkcolor=blue,          
    citecolor=magenta,        
    filecolor=magenta,      
    urlcolor=cyan           
}

\begin{document}
\title{On Polynomial Progressions Inside Sets of Large Dimension}

\author[B. Krause]{Ben Krause}
\address{Department of Mathematics, University of Bristol, BS8 1QU, UK}
\email{ben.krause@bristol.ac.uk}

\author[A. E. de Orellana]{Ana E. de Orellana}
\address{Department of Mathematics and Statistics, University of Jyv\"askyl\"a, FI-40014, Finland}
\email{ana.e.deorellana@jyu.fi}
\thanks{AEdO was supported by the European Research Council (ERC) under the European Union’s Horizon Europe research and innovation programme (grant agreement No 101087499).}

\newcommand{\ben}[1]{{\color{red}{#1}}}
\newcommand{\mariusz}[1]{{\color{blue}{#1}}}
\newcommand{\sarah}[1]{{\color{green}{#1}}}

\date{\today}

\begin{abstract}
In this note we connect Sobolev estimates in the context of polynomial averages e.g.\
\begin{align}\label{e:savings}
    \| \int_0^1 \prod_{k=1}^m f_k(x-t^k) \|_{1} \leq \text{Const} \cdot 2^{-\text{const} \cdot l} \prod_{i=1}^m \| f_k \|_m
\end{align}
whenever some $f_i$ vanishes on $\{ |\xi| \leq 2^l \}$ to the existence of polynomial progressions inside of sets of sufficiently large Hausdorff dimension, in analogy with work of Peluse in the discrete context. Our strongest results builds off work of Becker-Krause and is as follows: suppose that $\mathcal{P} = \{a_{d_1} t^{d_1}, a_{d_2} t^{d_2},\dots, a_{d_k} t^{d_k}\}$ vanish at the origin at different rates, and that $E \subset [0,1]$ has sufficiently large Hausdorff dimension,
\[ 1 - \text{const}(\mathcal{P}) < \text{dim}_H(E) < 1
\]
and Hausdorff content bounded away from zero, sufficiently large in terms of its dimension. Then $E$ contains a non-trivial polynomial progression of the form
\begin{align}
    \{ x , x - a_{d_1} t^{d_1}, x - a_{d_2} t^{d_2}, \dots, x - a_{d_k} t^{d_k} \} \subset E, \; \; \; t \neq 0.
\end{align}

Also, using the Fourier spectrum, we provide a short proof that whenever $E$ supports a measure with both positive Fourier dimension and a sufficiently large Frostman condition, it necessarily contains a non-trivial generalized three-term arithmetic progression of the form
\[ \{ x, x - \gamma_1 t, x- \gamma_2 t\} \subset E, \; \; \; \gamma_i \in \mathbb{Q},\  t \neq 0.\]
\end{abstract}

\maketitle


\section{Introduction}
In recent years, much work has been done in the context of detecting non-linear polynomial progressions inside of subsets of the integers of vanishing upper density; the strongest results in this direction are as follows \cite{PP,SW}.

\begin{theorem}[Peluse-Prendiville $k=3$; Shao-Wang $k \geq 4$]\label{t:disc}
    Suppose that $A \subset \{ 1,\dots,N\}$ with $N$ sufficiently large, and let $k \in \mathbb{N}$. Suppose that $P_1,\dots,P_{k-1}$ are fixed non-constant polynomials with distinct degrees and zero constant coefficients. Then there exist constants $0 < c_k < 1$ so that whenever
    \begin{align}
        |A| \geq N (\log N)^{-c_k},
    \end{align}
necessarily $A$ contains a $k$-term non-linear polynomial progression of the form
\begin{align}\label{e:discpoly} \{x, x-P_1(n),\dots, x-P_{k-1}(n) \}, \;\; \; x,n \in \{1,\dots,N\}.
\end{align}
\end{theorem}
In recent collaboration with Mousavi-Tao-Ter\"{a}v\"{a}inen \cite{K++}, the first-named author proved an analogous result when $n + 1$ in \eqref{e:discpoly} is prime.

A key ingredient in Theorem \ref{t:disc} is an inverse theorem, which, roughly speaking, says that, whenever $|f_i| \leq 1$ and the degrees of $\{P_1,\dots,P_m\}$ are $d_1 < \dots < d_m$, the only way that
\begin{align}\label{e:inv}
    \| \frac{1}{N} \sum_{n \leq N} \prod_{i=1}^m f_i(x-P_i(n)) \|_{\ell^1([1,N^{d_m}])} \geq \text{Const} \cdot \delta N^{d_m},
\end{align}
for $N \geq \delta^{- \text{Const}}$ sufficiently large 
is if $f_i$ ``look like" the plane wave $e( \beta_i x)$ where
\begin{align}
    \| q \beta_i \|_{\mathbb{T}} \leq \text{Const} \cdot \delta^{-\text{Const}} N^{-d_i}
\end{align}
for some $q \leq  \delta^{-\text{Const}}$. By arguing as in \cite{K+}, see \cite{Kosz} and \cite{Wan} for the strongest results in this direction, this implies the so-called \emph{arithmetic Sobolev estimate}
\begin{align}
    \| \frac{1}{N} \sum_{n \leq N} \prod_{i=1}^m f_i(x-P_i(n)) \|_{\ell^1} \leq \text{Const} \big(2^{- \text{const} \cdot l} + N^{-\text{const}} \big) \prod_{i=1}^m \| f_i \|_{\ell^m}
\end{align}
whenever some $f_i$ vanishes on
\begin{align}
    \bigcup_{q \leq 2^l} \{ \beta \in {\mathbb{T}} : \| q \beta \|_{\mathbb{T}} \leq 2^l N^{-d_i} \};
\end{align}
the inverse Theorem \eqref{e:inv} first appeared in deep work of Peluse \cite{P1}, and was established by a subtle application of PET-induction, see \cite{BL}. From the perspective of convergence, these estimates are crucial in proving pointwise convergence of polynomial ergodic averages. 

The goal of this paper is to investigate Euclidean analogues of Theorem \ref{t:disc}; this problem has a long history in the context of sets of positive upper density, namely $A \subset \mathbb{R}$ so that
\begin{align}
    \limsup_{|I| \to \infty \text{ an interval}} \frac{|A \cap I|}{|I|} > 0
\end{align}
dating back to work of Bourgain \cite{B}, with the strongest results due to the first-named author, in collaboration with Mirek, Peluse, Wright \cite{KMPW}, see also \cite{D+} for a similar perspective.

On the other hand, we focus on the issue of pattern detection inside of \textit{sparse} subsets of the unit interval. A line of inquiry initiated in important work of {\L}aba-Pramanik in \cite{LP}, with subsequent work due to \cite{HLP, FGP}. In all of these efforts, the key was for a set to support a measure with (large) \textit{Fourier dimension}, namely
\begin{align}
    \text{dim}_{F}\mu := \sup \{ s\in[0,1] : |\hat{\mu}(\xi)| \leq \text{Const}(1+|\xi|)^{-s/2} \}.
\end{align}

To have more nuanced control on the behavior of the Fourier decay of the measure, we will consider a family of dimensions called the Fourier spectrum \cite{Fra}. This allows us to weaken the assumptions of \cite{LP} of requiring a measure with Fourier dimension from $2/3$ to positive.

For $s\geq0$ and $\theta\in(0,1]$, define the $(s,\theta)$-energies by
\begin{equation*}
    \J_{s,\theta}(\mu) = \left(\int_{\mathbb{R}} |\widehat{\mu}(\xi)|^\frac{2}{\theta} |\xi|^{\frac{s}{\theta}-1}\,d\xi\right)^{\theta},
\end{equation*}
and for $\theta=0$,
\begin{equation*}
    \J_{s,\theta}(\mu) = \sup_{\xi\in\mathbb{R}}|\widehat{\mu}(\xi)|^{2}|\xi|^{s}.
\end{equation*}
Then the Fourier spectrum of $\mu$ at $\theta\in[0,1]$ is
\begin{equation*}
    \dim_{F}^\theta \mu = \sup\{ s\in[0,1] : \J_{s,\theta}(\mu)<\infty \}.
\end{equation*}
For $\theta=0$ the Fourier spectrum recovers the Fourier dimension $\dim_{F}\mu = \dim_F^{0}\mu$, and for $\theta=1$ the Sobolev dimension $\dim_{S}\mu = \dim_F^{1}\mu$. As a function of $\theta\in[0,1]$, the Fourier spectrum is non-decreasing, concave, and continuous in $(0,1]$. See \cite{Mat, W} for a review of the relationship between different notions of dimension and the Fourier transform.

In Theorem~\ref{t:hdim}, we remove the Fourier dimension constraint, and consider sets with large Hausdorff dimension; for preliminary results in this direction, see \cite{Kold}. To state our results, we introduce the following definition.

\begin{definition}
    Let $\mathbf{P}(\mathcal{P})$ denote the statement that for   \begin{align}
        \mathcal{P} = \{ P_1,P_2,\dots,P_m \}, \; \; \; P_k(t) := \sum_{n = e_k}^{d_k} a_n t^n
    \end{align}
    the following Sobolev estimate holds for some absolute $0 < \sigma_{\mathcal{P}} < C_{\mathcal{P}} < \infty$
 \begin{align}\label{e:sob}
    \| \int_0^1 \prod_{i=1}^m f_i(x-P_i(2^{-l}t)) \ dt \|_1 \leq C_{\mathcal{P}} N^{-\sigma_{\mathcal{P}}} \prod_{i=1}^m \|f_i \|_m
\end{align}
whenever some $f_i$ has Fourier transform that vanishes on
\[ \{ |\xi| \leq N \cdot 2^{l e_i} \}, \; \; \; l \geq C_{\mathcal{P}} \text{ sufficiently large}.\]
In the above notation, we will say $\mathcal{P}$ is \emph{relatively curved} if $\{ e_k \}$ are distinct, and we will refer to the constant $\sigma_{\mathcal{P}}$ as the \emph{Sobolev gain}.
\end{definition}
 In point of fact, using the improving nature of polynomial averaging operators, \eqref{e:sob} is equivalent up to constants to the above estimate with $L^{\infty}$ norms replacing $L^m$ norms, provided the integral on the left-hand side is taken over the unit interval; see e.g.\ \cite{CDR, K+} for details. As such, if $\mathbf{P}(\mathcal{P})$ holds, then necessarily $\mathbf{P}(\mathcal{P}')$ does as well, whenever $\mathcal{P}' \subset \mathcal{P}$. 

The current state of the art, for general polynomials, was established in deep work of Hu-Lie \cite{HL}, see \cite[Theorem 3.1]{HL} and Observation 1.2 (i); significant progress on the $m=2$ case was established in \cite{Li}, see also \cite{Lie0}.
\begin{theorem}[\cite{HL}]\label{t:tool}
    $\mathbf{P}(\mathcal{P})$ holds for all $|\mathcal{P}| \leq 3$.
\end{theorem}

Moreover, it was recently found by Becker and the first-named author in \cite{BK} that $\mathbf{P}(\mathcal{P})$ holds
monomial curves $\mathcal{P} = \{P_1, \dots, P_m\}$, $P_i(t) = a_it^{d_i}$, $1 \le d_1 < \dots < d_m$, of arbitrary length $m$. 
As a corollary of Theorem~\ref{t:tool}, by a Taylor expansion argument, see \cite[Remark 2]{Lie}, the following maximal estimate presents, which leads to the corresponding convergence result by an elementary density argument; for the strongest estimate of this type, see the recent paper \cite{GLZ}.

\begin{cor}
    The following maximal estimate holds whenever $\mathcal{P}$ is relatively curved with $n \leq 3$ polynomials:
    \begin{align}
        \| \sup_{k \geq 1} |B_k^{\mathcal{P}}(f_1,\dots,f_n)| \|_{1} \leq C_\mathcal{P} \prod_{i \leq n} \|f_i \|_n
    \end{align}
where
\begin{align}
    B_k^{\mathcal{P}}(f_1,\dots,f_n) := \int_0^{1} \prod_{i \leq n} f_i(x-P_i(2^{-k} t)) \ dt,
\end{align}
see \eqref{e:Bk} below.
\end{cor}

The aim of this paper is to connect $\mathbf{P}(\mathcal{P})$ to the detection of polynomial configurations inside sets of sufficiently large Hausdorff dimension -- but possibly zero Fourier dimension; see \cite[Proposition 1.3]{Zhu} for a similar statement and conclusion. 

In particular, our main result is as follows,

\begin{theorem}\label{t:hdim}
Let $\mathcal{P} = \{ P_1,\dots,P_m \}$ be a relatively curved collection of polynomials such that $\mathbf{P}(\mathcal{P})$ holds. Suppose that $A \subset [0,1]$ is a compact set with sufficiently large Hausdorff content
\begin{align}\label{e:cont}
    H^s_{\infty}(A) \geq a_{\mathcal{P}} > 0, \; \; \; s > 1 - c_{\mathcal{P}}.
\end{align}
Then, there exists an absolute constant $C_m$ so that whenever
\[ c_{\mathcal{P}} \leq C_m \sigma_{\mathcal{P}} a_{\mathcal{P}}^{10m} \]
is sufficiently small, there exist $x, t \in [0,1]$ so that $A$ contains the $m+1$ term polynomial progression
\[ \{x,x-P_1(t),x-P_2(t),\dots, x-P_m(t)\} \subset A, \; \; \; t \neq 0.\]
\end{theorem}

The proof of the theorem is by an energy pigeon-holing argument, and we begin by proving our result when $A$ has positive measure; this follows by specializing $f = \mathbf{1}_A$ in the Proposition \ref{p:SETUP} below. This proposition is a multi-linear analogue of Bourgain's work \cite{B} on non-linear Roth theorems and its subsequent extension \cite{D+}; with this result in hand, the result is proven by directly appealing to $\mathbf{P}(\mathcal{P})$.

In the course of completing this work we learned of a recent preprint \cite{Hong}, which uses different methods to address the periodic context under the assumption that the degrees of $\{P_1,\dots,P_m\}$ are distinct, adapting estimates from \cite{KMPW}. As is (essentially) observed there, the main result of \cite{Hong} follows from $\mathbf{P}(\mathcal{P})$, so in particular we recover the four-term arithmetic progression case of \cite{Hong} without the assumption on distinct degrees of the pertaining polynomials -- provided the collection is relatively curved.

The hypothesis of relative curvature is in general necessary to obtain the Sobolev savings of $\sigma_{\mathcal{P}}$, as a latent (approximate) modulation invariance presents when at least two polynomials vanish to the origin at the same degree. In the simplest such case, however, we are still able to detect non-trivial progressions, under the additional assumption that our set support a measure with both positive Fourier dimension and a sufficiently large Frostman condition, and our polynomials have rational coefficients.

In the following results, the dependence of $\varepsilon_{0}$ in the constants involved is essential, see \cite{Shm}. For this reason, we state our results in terms of finiteness of the $(s,\theta)$-energies. However, condition (1) can always be interpreted as a condition on the dimension of the measure, i.e. $\J_{s,\theta}(\mu) <\infty$ implies $\dim_{F}^\theta \mu \geq s$.

\begin{theorem}\label{t:main0}
  Let $\gamma_1,\gamma_2$ be non-zero rationals. Suppose $A\subset[0,1]$ is compact and supports a Frostman measure $\mu\in \text{Pr}(A)$ with
  \begin{enumerate}
        \item[(1)] $|\widehat{\mu}(\xi)|^2\leq C_1|\xi|^{-s}$ for some $s>0$; and
        \item[(2)] $\mu(B(x,r))\leq C_2r^{1-\varepsilon_0}$ for all $x\in\RR$, $r>0$, and some small enough 
        \[ \varepsilon_0 = \varepsilon_0(C_1,C_2,s,\gamma_1,\gamma_2).\]
    \end{enumerate}
Then $A$ contains a non-trivial 3AP of the form
\begin{align}
    \{ x,x-\gamma_1 t, x- \gamma_2 t\} \subset A, \; \; \; t \neq 0.
\end{align}
\end{theorem}
This estimate improves quantitatively on an analogous result of \cite{LP}, who required that $A$ to support a measure with Fourier dimension greater than $2/3$. 

Theorem~\ref{t:main0} is a consequence of the following more general statement in terms of the Fourier spectrum after setting $\theta=0$ and noting that if $\mu(B(x,r))\leq Cr^{\alpha}$, then $\J_{u,1}(\mu)<\infty$ for any $u<\alpha$.
\begin{theorem}\label{t:FS3APs}
   Let $\gamma_1,\gamma_2$ be non-zero rationals.  Suppose $A\subset[0,1]$ is compact and supports a Frostman measure $\mu\in \text{Pr}(A)$ with
  \begin{enumerate}
        \item[(1)] For some $\theta\in[0,1]$, $\J_{s,\theta}(\mu), \J_{u,1-\frac{\theta}{2}}(\mu) <\infty$ for $u,s>0$ satisfying $\frac{s}{2}+u>1$; and
        \item[(2)] $\mu(B(x,r))\leq Cr^{1-\varepsilon_0}$ for all $x\in\RR$, $r>0$, and some small enough 
        \[ \varepsilon_0 = \varepsilon_0(C,\J_{s,\theta}(\mu), \J_{u,1 - \frac{\theta}{2}}(\mu),\gamma_1,\gamma_2).\]
    \end{enumerate}
    Then $A$ contains a non-trivial 3AP of the form
\begin{align}
    \{ x,x-\gamma_1 t, x- \gamma_2 t\} \subset A, \; \; \; t \neq 0.
\end{align}
\end{theorem}

Optimizing condition (1) in Theorem~\ref{t:FS3APs} over $\theta\in[0,1]$ we obtain the following simpler statement, which pleasantly connects our numerology to that of \cite{LP}.
\begin{cor}
   Let $\gamma_1,\gamma_2$ be non-zero rationals. Suppose $A\subset[0,1]$ is compact and supports a Frostman measure $\mu\in\text{Pr}(A)$ with
  \begin{enumerate}
        \item[(1)] $\J_{s,\frac{2}{3}}(\mu) <\infty$ for some $s>\frac{2}{3}$; and
        \item[(2)] $\mu(B(x,r))\leq Cr^{1-\varepsilon_0}$ for all $x\in\RR$ $r>0$, and some small enough 
        \[ \varepsilon_0 = \varepsilon_0(C,\J_{s,\frac{2}{3}}(\mu),\gamma_1,\gamma_2).\]
    \end{enumerate}
  Then $A$ contains a non-trivial 3AP of the form
\begin{align}
    \{ x,x-\gamma_1 t, x- \gamma_2 t\} \subset A, \; \; \; t \neq 0.
\end{align}
\end{cor}

\bigskip

The structure of the paper is as follows: Theorem \ref{t:hdim} will be our first focus, which we address in \S 2-3, before addressing Theorem \ref{t:FS3APs} in \S 4.

\subsection{Acknowledgements}
We would like to thank Tuomas Orponen for carefully reading an earlier version, identifying the role of Hausdorff content, and patiently explaining things, as well as for providing the references \cite{Mat} and \cite{Zhu}. Thanks also go to Lars Becker for helping improve the exposition and overall clarity. We would also like to thank Jonathan Fraser, Alex Iosevich and Hamed Mousavi for helpful conversations, and Hamed in particular for his help determining the lower bound in Lemma \ref{l:lbd}.

\subsection{Notation}
Below, $\varphi$ will be mean-one non-negative Schwartz functions, and $\psi$ will denote Schwartz functions whose Fourier transforms are supported on an annulus away from the origin. We use the following notation for $L^1$-normalized dilations
\[ \phi_k(x) := 2^k \phi(2^k x), \]
and define
\begin{align}\label{e:Bk}
    B_k^{\mathcal{P}}(f_1,\dots,f_m)(x) := \int_0^1 \prod_{i=1}^m f_i(x-P_i(2^{-k} t)) \ dt.
\end{align}

We will make use of the modified Vinogradov notation. We use $X \lesssim Y$ or $Y \gtrsim X$ to denote
the estimate $X \leq CY$ for an absolute constant $C$ and $X, Y \geq 0.$  If we need $C$ to depend on a
parameter, we shall indicate this by subscripts, thus for instance $X \lesssim_{\mathcal{P}} Y$ denotes the estimate $X \leq C_{\mathcal{P}} Y$ for some $C_{\mathcal{P}}$ depending on $\mathcal{P}$. We use $X \approx Y$ as shorthand for $Y \lesssim X \lesssim Y$. We use the notation $X \ll Y$ or $Y \gg X$ to denote that the implicit constant in the $\lesssim$ notation is extremely large, and analogously $X \ll_{\mathcal{P}} Y$ and $Y \gg_{\mathcal{P}} X$.

We also make use of big-Oh and little-Oh notation: we let $O(Y)$  denote a quantity that is $\lesssim Y$ , and similarly
$O_p(Y )$ will denote a quantity that is $\lesssim_p Y$; we let $o_{t \to a}(Y)$
denote a quantity whose quotient with $Y$ tends to zero as $t \to a$ (possibly $\infty$), and
$o_{t \to a;p}(Y)$
denote a quantity whose quotient with $Y$ tends to zero as $t \to a$ at a rate depending on $p$.

\section{Energy Pigeon-holing}

The main result of this section is the following proposition; throughout, we will assume that $\mathbf{P}(\mathcal{P})$ holds.

\begin{proposition}\label{p:SETUP}
Suppose that $\mathcal{P}$ is relatively curved and that $0 \leq f \leq \mathbf{1}_{[0,1]}$ satisfies
\[ \int f \geq \epsilon > 0.\]
Suppose that $f_i = \varphi_{t_i}*f, \ 1 \leq i \leq m, \ 0 < t_i \leq \infty$,\footnote{We interpret $\varphi_{\infty} = \delta$, the point-mass at the origin.}
 and $0 \leq f_0 \leq \mathbf{1}_{[0,1]}$ has $\int f_0 \geq \epsilon$.
Then there exists some $0 \leq k \leq  \epsilon^{-10 m}$ so that
\[ \int f_0 \cdot B_k^{\mathcal{P}}(f_1,\dots,f_m) \gtrsim \epsilon^{m+1}.\]
\end{proposition}

By re-scaling the previous proposition, we arrive at the following corollary, which we will use below.
\begin{cor}\label{c:SETUP}
Suppose that $0 \leq f \leq M \cdot \mathbf{1}_{[0,1]}$, and that $\int f \geq 1$. Then
\[ \int f_0 \cdot B_0^{\mathcal{P}}(f_1,\dots,f_m) \gtrsim 2^{- M ^{10 m}} \]
whenever $f_i = f*\varphi_{t_i}$ as above and $0 \leq f_0 \leq M \cdot \mathbf{1}_{[0,1]}$ has integral $\geq 1$.
\end{cor}

The key ingredient in our proof of Proposition \ref{p:SETUP}, and thus Corollary \ref{c:SETUP}, is the following lemma, which we will iterate at many scales.

\begin{lemma}\label{l:TOOSMALL}
Let $k \geq 0$, and suppose that $0 \leq f_0,f_1,\dots,f_m \leq \mathbf{1}_{[0,1]}$ have $\int f_i \geq \epsilon$, but that
\[ \int f_0 \cdot B_k^{\mathcal{{P}}}(f_1,\dots,f_m) \ll \epsilon^{m+1}.\]


Then there exists some index $1 \leq i \leq n$ and some $|l| \lesssim \log(1/\epsilon)$ so that
\[ \| \psi_{ke_i+l}* f_i \|_m \gtrsim \epsilon^{1+2/m}.\]

\end{lemma}
\begin{proof}
Suppose that 
\begin{equation}\label{e:TOOSMALL}
\int f_0 \cdot B_k^{\mathcal{P}}(f_1,\dots,f_m) \ll \epsilon^{m+1}.
\end{equation}
With $s = C \cdot \log(1/\epsilon)$ for some sufficiently large $C$, we may expand
\begin{align} 
&\int f_0 \cdot B_k^{\mathcal{{P}}}(f_1,\dots,f_m)(x) = \int f_0 \cdot \prod_{i=1}^m \varphi_{k e_i-s}*f_i(x) \\
& \qquad + \sum_{|s_i| \leq s} \int f_0 \cdot B_k^{\mathcal{{P}}}(\psi_{s_1+ke_1}*f_1,\dots,\psi_{s_m+ke_m}*f_m)(x)  + O(\epsilon^{10m}), \end{align}
where we have used Lemma \ref{l:UncertP} below
to bound the low-frequency error terms  pointwise by an error of $\lesssim \epsilon^{10m}$, and \eqref{e:sob} to bound the high frequency terms.

Since \eqref{e:TOOSMALL} is small, while the low frequency terms are large by Lemma \ref{l:TECH}, by subtracting appropriately we must have a lower bound:
\[ \aligned 
\epsilon^{m+1} &\lesssim \sum_{|s_i| \leq s} \int |f_0| \cdot |B_k^{\mathcal{P}}(\psi_{s_1+ke_1}*f_1,\dots,\psi_{s_m + k e_m}*f_m)| \\
&\qquad \lesssim \sum_{|s_i| \leq s} \| B_k^{\mathcal{P}}(\psi_{s_1+k e_1}*f_1,\dots,\psi_{s_m + k e_m}*f_m)\|_1
\endaligned \]
from which the result follows from pigeon-holing and convexity.
\end{proof}

\begin{lemma}\label{l:UncertP}
Suppose $l \geq 0$, and that $\supp \hat{f_i}$ is supported in $\{ |\xi| \lesssim 2^{-l} 2^{k e_i}\}$.

Then
\[ 
B_k^{\mathcal{P}}(f_1,\dots, f_m) = f_i \cdot B_k^{\mathcal{P}_i}(f_1,\dots,f_{i-1}, f_{i+1},\dots,f_m) + \mathcal{E}_l,\]
where
\[ |\mathcal{E}_l| \lesssim 2^{-l} \prod_i \|f_i \|_{L^{\infty}}.\]
Here, $\mathcal{P}_i = \{ P_1,\dots,P_{i-1},P_{i+1},\dots,P_m\}$.
\end{lemma}
\begin{proof}
Let $\varphi$ be a Schwartz function with $\hat{\varphi} = 1$ on the support of $\hat{f_i}$, and vanishes in a neighborhood of it.

We set
\[ \mathcal{E}_l := B_k^{\mathcal{P}}(f_1,\dots,\varphi*f_i,\dots, f_m) - 
\varphi * f_i \cdot
B_k^{\mathcal{P}_i}(f_1,\dots,f_{i-1}, f_{i+1},\dots,f_m),\]
and use the mean-value theorem. More precisely,
\begin{align}
&B_k^{\mathcal{P}}(f_1,\dots,\varphi*f_i,\dots,f_m)(x) \\
& \qquad  = 
\int \left( 2^k \int_0^{2^{-k}} \prod_{j\neq i} f_j(x-P_j(t)) \cdot \varphi(x-P_i(t) - y)  \ dt \right) \ f_i(y) \ dy \\
& \qquad \qquad = 
\left( 2^k \int_0^{2^{-k}} \prod_{j\neq i} f_j(x-P_j(t)) \ dt \right) \cdot \int \varphi(x-y) f_i(y) \ dy \\
& \qquad \qquad \qquad + O\left(2^{-l} \cdot M_{\text{HL}} f_i \cdot \sup_{k} \,  B_k^{\mathcal{P}_i}(|f_1|, \dots, |f_m| ) \right), \end{align}
since
\[ \sup_{|u| \lesssim 2^{-ke_i}} |\varphi(x-u) - \varphi(x)| \lesssim 2^{-l} \cdot 2^{k e_i-l} \cdot (1 + 2^{k e_i-l}\cdot |x|)^{-100} \]
by the mean-value theorem. Above, $M_{\text{HL}}$ is the Hardy-Littlewood maximal function.
\end{proof}

As mentioned above, we make use of the following technical lemma.
\begin{lemma}\label{l:TECH}
For any $0 \leq f \leq \mathbf{1}_{[0,1]}$, any $k_i \geq 0$ and any $m$
\[ \int \prod_{i=0}^m \varphi_{k_i}*f \gtrsim \left( \int f\right)^{m+1}.\]
\end{lemma}
\begin{proof}
Following the approach of \cite{D+}, we minorize the convolution operators $\varphi_{k}*f$ with the conditional expectation operators
\[ \mathbb{E}_{k} f := \sum_{|I| = 2^{-k} \text{ dyadic}} \left( \frac{1}{|I|} \int_I f \right) \mathbf{1}_I,\]
at which point the result follows by induction.
\end{proof}

We now prove Proposition \ref{p:SETUP}.

\begin{proof}[Proof of Proposition \ref{p:SETUP}]
We proceed by contradiction, and in the notation of Proposition \ref{p:SETUP}, assume that for each $0 \leq k \leq K$ the following upper bound held
\[ \int f_0 \cdot B_k^{\mathcal{P}}(f_1,\dots,f_m) \ll \epsilon^{m+1}; \]
we will show that $K \lesssim \epsilon^{-10 m}$.

To do so, we use Lemma \ref{l:TOOSMALL} to extract an index $1 \leq i \leq m$, and a sparse subset of scales $X \subset \{ k \leq K \}$ separated by $\gg \log(1/\epsilon)$, of size
\[ |X| \gtrsim \frac{K}{m \cdot \log(1/\epsilon)},\]
with the following property:

For each $k \in X$ there exists some perturbation $|s_k| \lesssim \log(1/\epsilon)$ so that
\begin{align}\label{e:ENERGY} \| \psi_{s_k + k e_i} * f_i \|_m \gtrsim \epsilon^{1+2/m}. \end{align}
In particular, taking a $\ell^m$ sum of \eqref{e:ENERGY} yields the upper bound
\begin{align}\label{e:lmsum}
\epsilon^{5m} \cdot |X|  &\lesssim \sum_{k \in X} \| \psi_{s_k + k e_i} * f_i \|_m^m = \| \big( \sum_{k \in X} |\psi_{s_k + k e_i} * f_i |^m \big)^{1/m} \|_m^m \\
& \qquad \leq \| Sf_i \|_m^m  \lesssim_m \|f_i \|_m^m  \leq 1, 
\end{align}
where
\begin{align}\label{e:SFXN}
S f := (\sum_k |\psi_k* f|^2)^{1/2}
\end{align}
is the Littlewood-Paley square function, which is bounded on $L^p$ for $1 < p < \infty$. In particular, we have exhibited the desired upper bound, $K \lesssim \epsilon^{-10 \cdot m}.$
\end{proof}

With these preliminaries in mind, we now turn to the proof of Theorem \ref{t:hdim}

\section{The Proof of Theorem \ref{t:hdim}}
We begin by selecting a measure on $A$ so that we have the explicit bound,
\begin{align}\label{e:muLambda} \mu(I) \leq \Lambda \cdot |I|^{\beta}, \; \; \; I \cap A \neq \emptyset, \ 0 < \Lambda < \infty,\end{align}
where
\[ \beta < \dim_H(A) < 1 \]
is very close to $1$. 

By appropriately normalizing \cite[Theorem 8.8]{Mat}, we may bound
\begin{align}
    \Lambda \lesssim H_{\infty}^{\beta}(A)^{-1}, 
\end{align}
so by our assumption \eqref{e:cont}, we have an absolute bound
\[ \Lambda \lesssim a_{\mathcal{P}}^{-1} \lesssim_{\mathcal{P}} 1.\]

We track this dependence below, and emphasize that
\[ \| \varphi_{k}* \mu\|_\infty \lesssim \Lambda \cdot 2^{k(1-\beta)}\]
for each $k \geq 0$; the implicit constant is determined only by the Schwarz normalization on $\{ \varphi \}$.

\begin{proof}[The Proof of Theorem \ref{t:hdim}]
Seeking a contradiction, suppose that 
\[ A^m \cap \{ (x,x-P_1(t), \dots, x - P_m(t)) : 0 \leq x, t \leq 1 \} = \emptyset; \]
by compactness, we could find an integer $J$ sufficiently large so that the two sets are separated by
\begin{align}
    \gg 2^{-\sqrt{J}}.
\end{align}
So, with $\mu \in \text{Pr}(A)$ an appropriate Frostman measure, set 
\[ f := \varphi_{J}*\mu;\]
by a dyadic decomposition in physical space, it suffices to exhibit upper and lower bounds -- independent of $J$ -- for
\begin{equation}\label{e:UPDOWN}
\int f \cdot B_0^{\mathcal{P}}(f,\dots,f).
\end{equation}

To do so, with $1 \ll_{m, \Lambda} l$ a sufficiently large integer to be determined, decompose
\[ \aligned 
B_0^{\mathcal{P}}(f,\dots,f) &=  B_0^{\mathcal{P}}(\varphi_{l}*f,\dots,\varphi_{l}*f) +
\sum_{s> l} \left( \sum_{\max\{s_i\} =s } B_0^{\mathcal{P}}(\psi_{s_1}*f,\dots,\psi_{s_n}*f) \right) \\
& \qquad =: B_l + \sum_{s > l} B_{s}  . \endaligned \]
We will dictate that $\psi, \varphi$ have compact support in Fourier space.

Now, let $ \rho$ be a Schwartz function whose Fourier transform is one on the $m$-fold iterated sum set
\begin{align}\label{e:rho0}
W := \text{supp } \varphi + \dots + \text{supp } \varphi  \ \bigcup \ \text{supp } \psi + \dots + \text{supp } \psi  
\end{align}
and which vanishes in a neighborhood of this set:
\begin{align}\label{e:rho}
\mathbf{1}_{W} \leq \hat{\rho} \leq \mathbf{1}_{3W}.
\end{align}

Upon taking inner products, we may express \eqref{e:UPDOWN} as
\[  
\eqref{e:UPDOWN} = \langle \rho_{ l}*f, B_l \rangle +  \sum_{{s} > l} \langle \rho_{s}*f, B_{{s}} \rangle. \]
By applying the convexity bound
\[ \| \psi_l*f \|_m \lesssim \| \psi_l*f \|_{\infty}^{1-1/m} \lesssim 2^{l \cdot (1-\beta) \cdot (1-1/m)}, \]
we may sum  
\begin{align}
\sum_{{s} > l} | \langle \rho_{s}*f, B_{{s}} \rangle| &\lesssim \sum_{{s} > l} \| \rho_{s}*f \|_\infty \cdot \|B_{{s}} \|_1 \\
&  \lesssim_m \Lambda \sum_{{s} > l}   2^{s (1-\beta)} \cdot 2^{-\sigma_{\mathcal{P}} {s}} \cdot \sum_{\max\{s_i\} = s} \prod_{i=1}^m \| \psi_{{s_i}}*f\|_m \label{e:gain} \\
& \lesssim_m \Lambda \sum_{s > l}  2^{s (1-\beta)} \cdot 2^{-\sigma_{\mathcal{P}} {s}} \sum_{0 \leq s_1 \leq \dots \leq s_m \leq s} 2^{(1-\beta)(1-1/m) \sum_{1}^m s_i}\\
& \lesssim \Lambda \sum_{s > l} 2^{s (1-\beta)} \cdot 2^{-\sigma_{\mathcal{P}} {s}} \cdot 2^{(1-\beta)(m-1) s + o_m(s)} \\
& \lesssim \Lambda \sum_{s > l} 2^{-s ( \sigma_{\mathcal{P}} - (1-\beta) m - o_m(1))} \\
& \lesssim \Lambda 2^{-l ( \sigma_{\mathcal{P}} - (1-\beta) m - o_m(1))},
\end{align}
where the gain in \eqref{e:gain} follows from $\mathbf{P}(\mathcal{P})$. In particular, if we choose $\beta = \beta(\mathcal{P})$ sufficiently close to $1$, so that e.g.\
\[ m \cdot (1 -\beta) < \sigma_{\mathcal{P}}/3,\]
and ensure $l \gg_{\mathcal{P}} \log m$ is sufficiently large, obtain the asymptotic
\begin{align}\label{e:almost}
\eqref{e:UPDOWN} = \langle \rho_{l}*f,B_l \rangle + O\left(2^{-\sigma_{\mathcal{P}}/2 \cdot l} \right).
\end{align}

As far as convergence is concerned, an upper bound for the first term is given by 
\begin{equation}\label{e:muUP}
\| \rho_{l}*f\|_\infty \lesssim m \cdot \Lambda \cdot 2^{l(1-\beta)},
\end{equation}
which is  bounded independent of $J$. As for our lower bound, we have
\[ \langle \rho_{l}*f,B_l \rangle \gtrsim 2^{- \big(m \Lambda \cdot 2^{l(1-\beta)}\big)^{10 m}}, \]
by applying Corollary \ref{c:SETUP} and a dyadic decomposition in physical space. 

We now specialize
\[ l \approx \frac{1}{100m(1-\beta)} \]
to obtain the bound
\begin{align}\label{e:done}
    \eqref{e:UPDOWN} \gtrsim 2^{ - C_m \alpha_{\mathcal{P}}^{-10m }} - 2^{- \frac{C_m \sigma_{\mathcal{P}}}{1-\beta} } \gtrsim 2^{ - C_m \alpha_{\mathcal{P}}^{-10}}
\end{align}
provided
\[ \sigma_{\mathcal{P}}  \alpha_{\mathcal{P}}^{10m} \gg_m c_{\mathcal{P}} > {1-\beta}. \]
2
We finally remark that we may at last replace
\begin{align}
    B_0^{\mathcal{P}} \longrightarrow B_{0;\kappa}^{\mathcal{P}}
\end{align}
where
\begin{align}
    B_{0;\kappa}^{\mathcal{P}}(f_1,\dots,f_m)(x) := \int_\kappa^1 \prod_{i=1}^m f_i(x-P_i(t)) \ dt
\end{align}
for $0 < \kappa = \kappa(\beta,\mathcal{P}) \ll 1$ sufficiently small, without distorting the lower bound \eqref{e:done}, as the Sobolev savings persist for these truncated averages. This allows us to conclude that our polynomial progressions are non-trivial.
\end{proof}

We now turn to Roth's Theorem.
\section{The Proof of Theorem \ref{t:FS3APs}}
In this section, we prove that any probability measure on $[0,1]$ satisfying the following two hypotheses:

\begin{enumerate}
    \item For some $\theta\in[0,1]$ and $u,s>0$ satisfying $\frac{s}{2} + u >1$, $\J_{s,\theta}(\mu)<\infty$ and $\J_{u,1-\frac{\theta}{2}}(\mu)<\infty$; and
    \item $\mu(B(x,r))\leq Cr^{1-\varepsilon_0}$ for all $x\in\RR$ $r>0$, and some small enough $\varepsilon_0$;
    \end{enumerate}
necessarily contains the generalized arithmetic progression
\[ \{ x,x-\gamma_1t, x- \gamma_2 t \} \subset \text{supp } \mu, \; \; \; \gamma_i \in \mathbb{Q}, \ t \neq 0.\]

Notice that any Salem measure $\mu$ with dimension $> 2/3$ satisfies the first of these conditions, in analogy with \cite{LP}.

Throughout this section, we will regard $M \in \mathbb{N}$ as arbitrary but fixed, and will consider generalized 3APs of the form
\[ x, x-\gamma_1 t, x-\gamma_2 t, \; \; \; t \neq 0\]
where
\begin{align}
    \gamma_i \in \mathcal{Q}_M := \big( \bigcup_{N \leq M} \mathbb{Z}/N \big) \cap [-M,M].
\end{align} 

We begin with a standard partition of unity in Fourier space,
\begin{align}
\widehat{\varphi}(\xi) + \sum_{l \geq 1} \widehat{\psi}(\xi/2^l) =: \widehat{\varphi}(\xi) + \sum_{l \geq 1} \widehat{\psi_l}(\xi) = 1
\end{align}
where $\varphi$ and $\psi$ are Schwartz, with $\widehat{\psi}$ supported near $|\xi| \approx 1$. For measures $\mu$, we set
\begin{align}
    \mu_{\leq l} := \varphi_l*\mu, \; \; \; \mu_l := \psi_l*\mu.
\end{align}

Since $\widehat{\mu_{l}}$ is supported on $|\xi|\approx 2^l$, by definition of the $(s,\theta)$-energies,
\begin{equation}\label{e:L2mu}
    \|\widehat{\mu_{l}}\|_{p} = \left( \int_{|\xi|\approx 2^l} |\widehat{\mu_{l}}(\xi)|^p |\xi|^{\frac{sp}{2} -1}2^{-l\big( \frac{sp}{2}-1 \big)}\,d\xi \right)^{\frac{1}{p}} \lesssim 2^{-l( \frac{sp}{2}-1)\frac{1}{p}}\J_{s,\frac{2}{p}}(\mu)^{1/2}.
\end{equation}

The central object of study below will be the forms
\begin{align}
    \Lambda_{\Gamma}(\mu) := \int \mu(\frac{\gamma_2 x - \gamma_1t}{\gamma_2 - \gamma_1}) \ d\mu(x) d\mu(t);
\end{align}
more generally, define
\begin{align}
\Lambda_{\Gamma}\big(h,\mu,\mu\big) := \int h(\frac{\gamma_2 x - \gamma_1t}{\gamma_2 - \gamma_1}) \ d\mu(x) d\mu(t);
\end{align}
we first quickly verify that $\Lambda_{\Gamma}(\mu)$ converges, by showing
\begin{align}\label{e:Cauchy}
\sum_{l \geq 1} |\Lambda_{\Gamma}\big(\mu_l,\mu,\mu\big)| < \infty  
\end{align}
converges absolutely.
\begin{proof}[Proof \eqref{e:Cauchy}]
By Fourier inversion, we may express
\begin{align}
\Lambda_{\Gamma}\big(\mu_l,\mu,\mu\big) &= \int_{|\xi| \approx 2^l} \widehat{\mu_l}(\xi) \widehat{\mu}(\xi \frac{\gamma_2}{\gamma_1-\gamma_2}) \widehat{\mu}( \xi \frac{\gamma_1}{\gamma_2-\gamma_1}) \ d\xi \\
& \lesssim_M \| \widehat{\mu_l} \|_{2/\theta} \cdot \| \widehat{\mu} \cdot \mathbf{1}_{|\xi| \approx_M 2^l} \|_{4/(2-\theta)}^2 \\
&\lesssim \J_{s,\theta}(\mu)^{1/2}\J_{u,1-\frac{\theta}{2}}(\mu) 2^{-l\big( \frac{s}{2}+u-1 \big)}
\end{align}
by hypothesis $(1)$, which yields \eqref{e:Cauchy} by summing a geometric series.
\end{proof}

We will prove Roth's Theorem for generalized $3APs$ holds, i.e.\ the support of $\mu$ contains a non-trivial progression
\[ x, x- \gamma_1 t, x- \gamma_2 t, \; \; \; t \neq 0, \ \gamma_i \in \mathcal{Q}_M\]
by proving that $\Lambda_{\Gamma}(\mu)>0$, and that the form gives zero weight to trivial $3APs$ (i.e.\ when $t = 0$ in the above).

Specifically, we will establish the following two propositions in turn:

\begin{proposition}\label{p:convform}
    With $\mu$ satisfying the assumptions $(1),(2)$ as above with $\varepsilon_{0}$ sufficiently close to $0$,
    \[ \Lambda_{\Gamma}(\mu) \gtrsim_M 1.\]
\end{proposition}

\begin{proposition}\label{p:diag}
    For each $\mu$ satisfying the assumptions $(1),(2)$ as above, whenever $\varepsilon_{0} < 1/2$ and $\Lambda_{\Gamma}(\mu) > 0$, $\mu$ necessarily contains a non-trivial progression in its support
    \[ \{x ,x - \gamma_1 t, x- \gamma_2 t \} \subset \text{supp}( \mu), \; \; \; t \neq 0. \]
\end{proposition}

We begin with Proposition \ref{p:convform}, which follows from transferring to the discrete situation. Specifically, we import the following quantitative bound from the integer setting.

\begin{lemma}\label{l:lbd}
    For every $0 \leq f \leq 1$ with $\int_0^1 f \geq \delta >0$, there exists an absolute constant, $c_M(\delta) > 0$, so that
    \begin{align}
        \int_{[0,1]^{2}} f(x) \prod_{i=1}^2 f(x-\gamma_i t) \ dx dt \geq c_{M}(\delta).
        \end{align}
In fact, one may choose
    \[ c_M(\delta) \gtrsim_M \exp(-\delta^{-1+c})
    > 0, \; \; \; 0 < c \ll 1.\]
\end{lemma}
\begin{proof}[Sketch]
By density, the above lower bound is realized by continuous functions (up to $\epsilon$ losses); since $\gamma_1,\gamma_2 \in \mathbb{Q}$, we may transfer to the discrete setting by a standard Riemann summation argument to export the quantitative bound. Indeed, applying \cite[Theorem 1.1]{P22} with $G=\mathbb{Z}_{C N}, \ C = C(M) \gg 1$, whenever $A \subset \{1,\dots,N\}$, has relative density
\[ \frac{|A|}{N} = \delta \gtrsim (\log N)^{-1-c}, \textup{ for some absolute constant } c>0, \]
necessarily $A$ contains a three-term arithmetic progression in $\mathbb{Z}_{CN}$, and since $A \subset \{1,\dots,N\} =: [N]$, there are no ``overlap" issues, provided $C$ is sufficiently large. This implies that 
\begin{align}
    N_3(\delta) &:= \min \big\{ N : \text{whenever $A \subset [N]$ has density $\geq \delta$,}\\
    & \qquad \qquad \qquad A \text{ contains a generalized 3AP} \big\}
\end{align} 
satisfies
\[ N_3(\delta) \lesssim \exp(\delta^{-1+c}),\]  
so adapting the proof in \cite[Theorem 18]{Shkredov} implies that 
\begin{align}
    c_M(\delta) \gtrsim_M \delta^{O(1)}N_{3}(\delta/2)^{-O(1)} \gtrsim \exp(-\delta^{-1+c/2}).
\end{align}
\end{proof}

\begin{proof}[Proof of Proposition \ref{p:convform}]
Abbreviate $\mu_0 := \mu_{\leq 0}$, and decompose
\begin{align}
    \mu = \mu_0 + \sum_{l \geq 1} \mu_l;
\end{align}
define
\begin{align}
    \Lambda_l(\mu) := \sum_{\max\{l_i\} = l} \Lambda_{\Gamma}\big(\mu_{l_1},\mu_{l_2}, \mu_{l_3} \big).
\end{align}
Note that by Fourier analysis   
\begin{align}
    \Lambda_l(\mu) &\leq \sum_{i,j = O_M(1)} \int_{|\xi| \approx 2^l} |\widehat{\mu_l}(\xi)| |\widehat{\mu_{l+i}}(\xi \frac{\gamma_2}{\gamma_1-\gamma_2})| |\widehat{\mu_{l+j}}(\xi \frac{\gamma_1}{\gamma_2-\gamma_1})| \ d\xi \\
    & \lesssim \J_{s,\theta}(\mu)^{1/2} \J_{u,1-\frac{\theta}{2}}(\mu) 2^{-l( \frac{s}{2}+u-1 )}.
\end{align}
So, with $l_0$ a free parameter, decompose
\begin{align}
    \Lambda_{\Gamma}(\mu) = \| \mu_{\leq l_0} \|_{{\infty}}^3 \Lambda_{\Gamma}\big( \frac{\mu_{\leq l_0}}{\| \mu_{\leq l_0} \|_{{\infty}}}\big) + \sum_{l > l_0} \Lambda_l(\mu).
\end{align}
We can bound the first term below by
\[ c_M(2^{-\varepsilon_{0}l_0}),\]
so
\begin{align}
    \Lambda_{\Gamma}(\mu) &\geq c_M(2^{-\varepsilon_{0} l_0}) - \J_{s,\theta}(\mu)^{1/2}\J_{u,1-\frac{\theta}{2}}(\mu) 2^{-l_{0}( \frac{s}{2}+u-1 )}\\
    & \geq \Omega_M(\exp( - \exp( \varepsilon_{0}l_0) )) - O(2^{-l_{0}( \frac{s}{2}+u-1 )}).
\end{align}
If we choose $l_0 = O_{M,s,u}(1)$ to be large but fixed, then whenever (say)
\begin{equation*}
    \varepsilon_{0}\leq \frac{1}{10 l_{0}}
\end{equation*}
we may bound
\begin{align}
    \Lambda_{\Gamma}(\mu)  \gtrsim_M 1,
\end{align}
uniformly in $\varepsilon_{0}$ in the above range.
\end{proof}

By Proposition \ref{p:convform}, we see that 
\[ \text{supp}(\mu)^3 \cap \{ (x,x-\gamma_1 t, x- \gamma_2 t) : 0 \leq x,t \leq 1 \} \neq \emptyset, \]
as otherwise, we could find $l_0$ so that the above sets would be separated by $\geq 2^{-\sqrt{l_0}}$. With $l \gg l_0$ a free parameter to be determined later, we could then decompose
\begin{align}
    1 \lesssim_M \Lambda_\Gamma(\mu) = \Lambda_\Gamma(\mu_{\leq l}) + \sum_{l' >l} \Lambda_l(\mu) = \Lambda_\Gamma(\mu_{ \leq l}) + O(2^{-\epsilon l});
\end{align}
since the first term on the right is $O(2^{-\sqrt{l}})$ (say) by a dyadic decomposition in physical space, we derive a contradiction upon sending $l \to \infty$.

We now complete the proof of Roth's Theorem by establishing Proposition \ref{p:diag}.

\begin{proof}[Proof of Proposition \ref{p:diag}]
With $\chi \geq 0$ a rapidly decaying Schwartz function with compact Fourier support, and
\begin{align}
    \chi^\delta(t) := \chi(\delta^{-1}(t-1)),
\end{align}
it suffices to show that
\begin{align}
    \Lambda_{\Gamma}\big(\chi^\delta \mu,\mu,\mu\big) = o_{\delta \to 0;\mu}(1),
\end{align}
i.e.\ that $\Lambda_{\Gamma}$ gives no mass to the diagonal. To see this, let $\delta>0$ be arbitrary, and telescope
\begin{align} 
\Lambda_{\Gamma}\big(\chi^\delta  \mu,\mu,\mu\big) = \Lambda_{\Gamma}\big(\chi^\delta \mu_{\leq l_0},\mu,\mu\big) + \sum_{l > l_0} \Lambda_{\Gamma}\big(\chi^\delta \mu_l,\mu,\mu\big).
\end{align}
We choose
\[ 2^{-l_0} := \delta^{\frac{1-\varepsilon_{0}-\kappa}{\varepsilon_{0}}}, \; \; \; 0 < \kappa < 1-2\varepsilon_{0}\]
so that the first term contributes
\[ O(2^{-\varepsilon_{0} l_0} \delta^{1-\varepsilon_{0}}) = O(\delta^{\kappa/2}) \]
by a dyadic decomposition in physical space.
Thus, the first term is $o_{\delta \to 0}(1)$; note that $2^{-l_0} \ll \delta$. To address the error terms, we prove geometric decay in $l > l_0$; specifically, we bound 
\begin{align}
    |\Lambda_{\Gamma}\big(\chi^\delta \mu_l,\mu,\mu \big)| &\leq \int_{|\xi| \approx 2^l} |\widehat{\mu}( \xi \frac{\gamma_2}{\gamma_1-\gamma_2})| |\widehat{\mu}(\xi \frac{\gamma_1}{\gamma_2-\gamma_1})| \left( \int_{|\zeta| \approx 2^l} \delta |\widehat{\chi}(\delta(\xi - \zeta))| |\widehat{\mu_l}(\zeta)| \, d\zeta \right) \, d\xi \\
    &\lesssim \| \widehat{\mu} \cdot \mathbf{1}_{|\xi| \approx_M 2^l} \|_{4/(2-\theta)}^2 \|\widehat{\chi}\|_{1} \| \widehat{\mu_{l}}\|_{2/\theta} \\
    &  \lesssim 2^{-l(\frac{s}{2}+u-1)}
\end{align}
where we can restrict the range of $\xi$ using the compact support assumption of $\widehat{\chi}$. 
Summing the geometric series yields a bound
\begin{align}
\Lambda_{\Gamma}\big(\chi^\delta \mu,\mu,\mu \big) = O(\delta^{\epsilon})
\end{align}
for some $\epsilon = \epsilon(s,u,\varepsilon_{0},M) > 0$.

The proof is complete.

\end{proof}

\typeout{get arXiv to do 4 passes: Label(s) may have changed. Rerun}

\end{document}